\documentclass[aps,prd,groupedaddress,amsmath,amssymb,amsfonts,nofootinbib,preprint,hidelinks]{revtex4-1}
\usepackage{amssymb}
\usepackage{amsmath, amssymb, ,amsfonts, verbatim, url, mathrsfs}
\usepackage{lineno}
\usepackage{mathrsfs}
\usepackage{mathtools}
\usepackage{verbatim}
\usepackage{amsthm}
\usepackage{framed}
\usepackage{wasysym}
\usepackage{upgreek}
\usepackage{color}
\usepackage[dvipsnames]{xcolor}
\usepackage{tensor}
\usepackage{accents}
\usepackage{dsfont}
\usepackage{hyperref}
\usepackage{enumerate}
\usepackage[normalem]{ulem}
\usepackage{longtable}
\usepackage{mathtools}
\usepackage{stmaryrd}

\numberwithin{equation}{section}

\newtheorem{theorem}{Theorem}[section]

\theoremstyle{definition}
\newtheorem{definition}{Definition}[section]
\newtheorem{remark}{Remark}[section]

\newcommand{\hta}{\tilde{\tau}}
\newcommand{\vertiiii}[1]{{\left\vert\kern-0.25ex\left\vert\kern-0.25ex\left\vert\kern-0.25ex\left\vert #1 \right\vert\kern-0.25ex\right\vert\kern-0.25ex\right\vert\kern-0.25ex\right\vert}}
\newcommand{\vertiii}[1]{{\left\vert\kern-0.25ex\left\vert\kern-0.25ex\left\vert #1 \right\vert\kern-0.25ex\right\vert\kern-0.25ex\right\vert}}

\newcommand{\Rbb}{\mathbb{R}}
\newcommand{\Zbb}{\mathbb{Z}}
\newcommand{\Tbb}{\mathbb{T}}
\newcommand{\Pbb}{\mathbb{P}}

\newcommand{\del}[1]{{\partial_{#1}}}

\newcommand{\AND}{{\quad\text{and}\quad}}

\newcommand{\la}{\langle}
\newcommand{\ra}{\rangle}

\newcommand{\p}[1]{
	\begin{pmatrix}
		#1
	\end{pmatrix}
}

\newcommand{\be}{\begin{equation}}
	\newcommand{\ee}{\end{equation}}

\allowdisplaybreaks

\begin{document}

	\title[Slightly nonlinear  Jeans instability]{Rigorous proof of the slightly nonlinear Jeans instability in the expanding Newtonian universe}

	\author{Chao Liu}
	\email{chao.liu.math@foxmail.com}
	\affiliation{Center for Mathematical Sciences and School of Mathematics and Statistics, Huazhong University of Science and Technology, Wuhan 430074, Hubei Province, China.}
	
	\author{Yiqing Shi}
	\email{yiqingshi@hust.edu.cn}
	\affiliation{School of Mathematics and Statistics, Huazhong University of Science and Technology, Wuhan 430074, Hubei Province, China.}

	\begin{abstract}
		Due to the nonlinearity of the Euler--Poisson equations, it is possible that the nonlinear Jeans instability may lead to a faster density growing rate than the rate in the standard theory of linearized Jeans instability, which \textit{motivates} us to study the nonlinear Jeans instability.  The \textit{aim} of this article is to develop a method proving the Jeans instability for slightly nonlinear Euler--Poisson equations in the expanding Newtonian universe. The standard proofs of the Jeans instability rely on the Fourier analysis. However, it is difficult to generalize Fourier method to a nonlinear setting, and thus there is no result in the nonlinear analysis of Jeans instability.
		We firstly develop a non-Fourier-based method to reprove the linearized Jeans instability in the expanding Newtonian universe. Secondly, we generalize this idea to a slightly nonlinear case. This method relies on the Cauchy problem of the Fuchsian system due to the recent developments of this system in mathematics.  The fully nonlinear Jeans instability for the Euler--Poisson and Einstein--Euler equations are in progress.
	\end{abstract}

	\maketitle
	
	\setcounter{tocdepth}{2}
	
	\pagenumbering{roman} \pagenumbering{arabic}

	\section{Introduction}

	Observations of the universe suggest that the small inhomogeneities in the density distribution  may lead to a well-developed nonlinear structure. A useful and widely used theoretical tool to investigate formations of protostars and the nonlinear structures, in astrophysics, is the famous Jeans instability, while the Jeans instability is established based on the \textit{linearized} Euler--Poisson system for the perturbed density, pressure, velocities, and the potential, and the standard proofs of the Jeans instability rely on the Fourier analysis. From a mathematical point of view, Fourier analysis can not be applied to the fully nonlinear equations. As a consequence,  there is no result in the nonlinear analysis of Jeans instability, which has been pointed out by Rendall in \cite{Rendall2002} since $2002$. In addition, we emphasize that the famous Jeans' criterion is a pointwise inequality for every time rather than a data of the Cauchy problem of the Euler--Poisson system,  since the Jeans length is time dependent for an expanding Newtonian universe (see \cite{Bonnor1957,ViatcehslavMukhanov2013} and it characterizes the competing of the pressure and gravity at every moment).
	
	In proofs of the linear Jeans instability by Fourier analysis, the smallness of the perturbed density is crucial, otherwise, the large density variations will significantly destroy the linearity of the Euler--Poisson equations. Therefore, the linear Jeans instability can only be valid for a small time region. We intend to develop a method which can lead to the fully nonlinear Jeans instability by giving the initial data, and it is expected to obtain a much faster growing rate of the nonlinear Jeans instability. We emphasize that our method on the nonlinear problem is not solving the approximation equations order by order\footnote{Usually, a solution to an approximation equation may not approximate to the solution of the original exact equation, although we always use this assumption to build the approximation methods. That is, the approximations in the equation level do not mean those in the solution level. }  and our proof is independent of the original classical Jeans' method. We want to directly solve the nonlinear equation, by using certain modern mathematical techniques, which will give a  precise estimate on the real solution to the original Euler--Poisson system rather than solutions of the approximation of the Euler--Poisson system. This article can not completely achieve this goal, but makes a first step to this aim and gives an estimate on the solution of the \textit{slightly} nonlinear Jeans equation [see \eqref{eq 18b}]. We intend to generalize this idea to the fully nonlinear case, i.e., proving the existence of the solution and estimating the exact solution to the nonlinear Euler--Poisson equations, in the near future.  
	
	In order to develop the \textit{nonlinear analysis} of Jeans instability for the Cauchy problems of the Euler--Poisson and Einstein--Euler systems, it is crucial to build a method which is independent of the Fourier analysis even for the linearized cases. The \textit{key motivation} of this article is to develop a method which is more adapted to the \textit{nonlinear analysis} of the Jeans instability of the Cauchy problems. That is, we aim to explore some significant structures of the equations which can also capture the nonlinear structures of the Jeans instability. In the progressing article, we will develop this method to the fully nonlinear Jeans instability in a Cauchy problem without the pointwise Jeans criterion.

	Jeans instability can be obtained by various methods (see, e.g.,  \cite{Zeldovich1971, Ward-Thompson2011,Capozziello2012,ViatcehslavMukhanov2013,Arbuzova2014}). The governing equations of Jeans instability is the Euler--Poisson system for isentropic fluids which characterizes   the fluid-filled universe with Newtonian self-gravity,
	\begin{align}
		\left\{
		\begin{aligned}
			&\frac{\partial\rho}{\partial t}+\nabla(\rho\boldsymbol v) = 0 ,\\
			&\frac{\partial\boldsymbol v}{\partial t}+(\boldsymbol v\cdot\nabla)\boldsymbol v+\frac{\nabla p}{\rho}+\nabla\phi = 0 ,\\
			&\Delta\phi =  4\pi G\rho,  \label{eq 3}	
		\end{aligned}
		\right.
	\end{align}
	where $\rho(t,\boldsymbol x)$, $\boldsymbol v(t,\boldsymbol x)$, $p(\rho)$, and $\phi$ are the energy density, 3-velocities, pressure of the fluids, and gravitational potential, respectively.

	The long time stability problem of Friedmann-Lema{\^{i}}tre-Robertson-Walker-like (i.e., FLRW-like for short)  universe has been widely investigated for linear and even for nonlinear cases with the small data perturbations in mathematical literature; see, for example, \cite{Brauer1994} for the nonlinear Newtonian cosmology and \cite{Ringstroem2008,RodnianskiSpeck:2013,Hadzic2015,Luebbe2013,Oliynyk2016a,Liu2018,Liu2018b,Liu2018a} for the nonlinear general relativity. We point out that \cite{Liu2018a} implies, for polytropic gas $p(\rho)=\kappa\rho^{\frac{n+1}{n}}$ $[n\in(1,3)]$, if the initial background density is smaller than certain bounds, it can conclude the global nonlinear stability of the FLRW-like universe for time $t\in[0,\infty)$. This leaves rooms and possibilities for the nonlinear Jeans instability. We remark that the nonlinear analysis \cite{Liu2021d,Liu2021} of the Euler--Poisson system implies that phenomena of mass accretions are not only due to the gravitational collapse, but can also be due to the effects of the free-falling boundary for irregular shaped molecular clouds.

	In order to avoid the famous \textit{Jeans swindle}\footnote{That is, the zero order approximations (i.e., the homogeneous background solutions) are not a solution of the zero order Euler--Poisson equations. See \cite{Falco_2013} for details. } of the original Jeans instability (see, e.g., \cite[\S$6.2$]{ViatcehslavMukhanov2013}, and \cite[\S$10.2$]{Zeldovich1971}), we consider an expanding homogeneous and isotropic Newtonian universe according to the assumptions\footnote{Note that alternative assumptions and formulations of the expanding universe can be found in \cite{Bonnor1957}. } in  \cite[\S$6.3$]{ViatcehslavMukhanov2013}, the background velocities, we assume, obey the \textit{Hubble law}:
	\begin{equation}
		\rho=\rho_{0}(t),\quad
		\boldsymbol v=\boldsymbol v_{0}=H(t) \boldsymbol x.\label{eq 4}
	\end{equation}
	Then, substituting\footnote{Take the divergence of the momentum equation in the Euler equation and use the Poisson equation to obtain the second equation in \eqref{eq 5}. } \eqref{eq 4} into \eqref{eq 3}, we arrive at the conservation of the total mass and Friedmann equation (see \cite[\S$6.3$]{ViatcehslavMukhanov2013} for details),
	\begin{align}\left\{
		\begin{aligned}
			&\dot{\rho_0}+3H\rho_{0}=0,\\
			&\dot{H}+H^{2}=-\frac{4\pi G}{3}\rho_0.\label{eq 5}
		\end{aligned}
		\right.
	\end{align}	
	First note that there is an exact solution to this homogeneous and isotropic Euler--Poisson \eqref{eq 5},
	\begin{align}
		\left\{
		\begin{aligned}
			&\rho_0(t)= \frac{1}{6\pi Gt^2},\\
			&H(t)=\frac{2}{3t}. \label{eq 8}
		\end{aligned}
		\right.
	\end{align}
 Then
	\begin{equation}\label{e:exsol}
		\rho_0(t)= \frac{1}{6\pi Gt^2}, \quad p_0(t)=f(\rho_0(t)), \quad \boldsymbol{v_0}(t,\boldsymbol{x})=\frac{2}{3t}\boldsymbol{x}, \AND \phi_0(t,\boldsymbol{x})=\frac{2}{3}\pi G\rho_0 |\boldsymbol{x}|^2,
	\end{equation}
where $f$ is a smooth, positive and increasing function, are the exact homogeneous solution\footnote{Also see \cite{Zeldovich1971} for details. } to \eqref{eq 3} in $(t,\boldsymbol{x})\in[1,\infty)\times \Rbb^3$.

	In order to analyze the behaviors of the perturbed variables deviating from the background solution \eqref{eq 8}, let us first decompose the variables $(\rho,\boldsymbol v,p,\phi)$ to the exact background solution \eqref{eq 8} and the perturbed parts as the following,
	\begin{equation}
		\rho =\rho_{0}+\delta\rho,\quad
		\boldsymbol v=\boldsymbol v_0 +\delta\boldsymbol v,\quad \phi=\phi_{0}+\delta\phi, \AND
		p =p_{0}+\delta p=p_{0}+c_s^2 \delta\rho,\label{eq 6}
	\end{equation}
	where $c_s^2=\mathrm{d}p/\mathrm{d}\rho$ is the square of the speed of sound.

	In this article, we first reprove, by a non-Fourier-based proof, the linearized Jeans instability for the expanding Newtonian universe, and improve the standard results.  Thus, we first recall,
	by introducing the fractional amplitude of the density perturbations $\varrho\equiv\delta\rho/\rho_0$, that $\varrho $ satisfies the linearized equation (see \cite[\S$6.3$]{ViatcehslavMukhanov2013}, and we call it the \textit{Jeans equation} for short in this article),
	\begin{equation}
		\ddot{\varrho}+2H\dot{\varrho}-\frac{c_s^2}{a^2}\Delta\varrho-4\pi\mathnormal{G}\rho_0\varrho=0 ,   \label{eq 7}
	\end{equation}
	where we note $\Delta$ and $\nabla$ now are the derivatives with respect to the Lagrangian coordinates $\boldsymbol q$ defined by $ \boldsymbol{x}=a(t)\boldsymbol q $, where\footnote{In fact, $a(t)=a(1) t^{\frac{2}{3}}=  t^{\frac{2}{3}}$ provided $a(1)=1$, since by the Hubble law \eqref{eq 4} and the Lagrangian coordinates $\boldsymbol{x}=a(t)\boldsymbol q$, we obtain  $H(t):=\frac{\dot{a}(t)}{a(t)}$. Then by \eqref{eq 8}, we can solve $a(t)$. } $a(1):=1$, and the time derivatives are obtained at constant $\boldsymbol q$ (i.e., the material derivatives).
	Let us briefly recall the derivations of this equation \eqref{eq 7} and  we refer readers to \cite{ViatcehslavMukhanov2013}, Sec. 6.3, for detailed derivations and equations. Substituting the decomposition \eqref{eq 6} into the Euler--Poisson system \eqref{eq 3}, with the help of \eqref{eq 5}, we are able to derive a set of linearized equations. Then using the Lagrangian coordinates $\boldsymbol q$, the linearized Euler--Poisson system of the perturbed variables can be further simplified by noting $\del{t}|_{\boldsymbol{x}}=\del{t}|_{\boldsymbol q}-\boldsymbol{v}_0\cdot \nabla_{\boldsymbol{x}}$ and $\nabla_{\boldsymbol{x}}=a^{-1}\nabla_{\boldsymbol{q}}$. At the end, taking the divergence of the above linearized momentum conservation, and using the continuity equation and the Poisson equation to replace $\nabla\cdot\delta\boldsymbol{v}$ and $\Delta\delta\phi$, respectively, we are able to conclude \eqref{eq 7}.

	For simplicity and in order to emphasize the main structure of the equations, in this article, we take a specific equation of state
	\begin{equation}\label{e:eos}
		p(\rho)=\kappa\rho^\gamma(\gamma>1)
	\end{equation}
	(i.e., the polytropic gas) to proceed the derivations. In this case, 	we further linearize \eqref{eq 7} by noting $c_s^2=\frac{\mathrm{d}p}{\mathrm{d}\rho}=\gamma\kappa\rho_0^{\gamma-1}(1+\varrho)^{\gamma-1}$ and $a(t)=t^{\frac{2}{3}}$, then \eqref{eq 7} becomes
	\begin{equation}
		\ddot{\varrho}+ \frac{4}{3t}\dot{\varrho}-\tilde{\kappa}t^{-2\gamma+\frac{2}{3}}\Delta\varrho- \frac{2}{3 t^2}\varrho=0\label{eq 9},\\
	\end{equation}
	where $\tilde{\kappa}=\gamma\kappa(\frac{1}{6\pi G})^{\gamma-1}$.

The second main topic of this article is the Jeans instability with a slightly nonlinear term in above linearized Jeans Eq. \eqref{eq 9} of $\varrho$. More precisely, this nonlinear Jeans equation is given by
	\begin{equation}
		\ddot{\varrho}+\frac{4}{3t} \dot{\varrho}-\tilde{\kappa}t^{-2\gamma+\frac{2}{3}}\Delta\varrho-\frac{2}{3t^2}\varrho=(\gamma-1)\tilde{\kappa}t^{-2\gamma+\frac{2}{3}}\frac{D^i\varrho D_i\varrho}{1+\varrho}.  \label{eq 18b}
	\end{equation}
In fact, we do not add the nonlinear term in the right hand of Eq. \eqref{eq 18b} arbitrarily, since this nonlinear term appears in the fully nonlinear Jeans equations (derived by the fully nonlinear Euler--Poisson system).

	The \textit{key tool} of this article is the Cauchy problem of the Fuchsian systems which was first established in \cite{Oliynyk2016a} in a nonlinear fashion and then has been investigated by a series works, for examples,  Oliynyk, Beyer, Olvera-Santamar{\'{\i}}a, and the first author of this paper in \cite{Liu2018,Liu2018b,Beyer2020,LW2021a}. We point out that the Fuchsian system in these works are the more general quasilinear ones, but in this article we only present the simplest semilinear formulations in Appendix \ref{s:Fuc}. We \textit{emphasize} that although the Fuchsian method may not be necessary in the linear case, it will be very promising for proof of the nonlinear Jeans instability, and we present a simple nonlinear case in this article.

    In this article, since we only focus on the perturbations of density by Eqs. \eqref{eq 9} and \eqref{eq 18b}, and  these density equations are independent of the velocity $\boldsymbol{v}$ and Newtonian potential $\phi$,  for simplicity, we solve these equations by restricting ourselves to the region $(t, \boldsymbol{q})\in M:=[1,\infty)\times \Tbb^3$, where $\Tbb^3:=S^1 \times S^1 \times S^1$ is a standard torus with period $1$. We point out that these methods can be generated to the region\footnote{In fact, for fully nonlinear case, the equation of density can not be fully decoupled, and the background potential $\phi_0$ given by \eqref{e:exsol} can not be defined on $\Tbb^3$. However, an alternative model by considering a Poisson equation with positive cosmological constant will leads to a uniform $\phi_0$ which can be defined on $\Tbb^3$. We leave these details to later works. } $(t, \boldsymbol{q})\in M:=[1,\infty)\times \Rbb^3$ with minor mathematical modifications.

In Sec. \ref{s:case2}, we rigorously reprove, for $\gamma=4/3$, the linearized Jeans instability in an expanding Newtonian universe by a non-Fourier-based method. In this method, a variation of Jeans criterion appears [given in \eqref{e:jeanscri}].  Then in the subsequent Sec. \ref{s:case1}, inspired by the method of Sec. \ref{s:case2}, we develop the Fuchsian formulation for the Eq. \eqref{eq 9} and estimate the behavior of the solution without a detailed solving of this equation. This section presents the main idea of this article and it can be applied to the nonlinear case. In the end, in Sec. \ref{s:case3}, we generalize the Fuchsian method in Sec. \ref{s:case1} to prove the slightly nonlinear Jeans instability for the Eq. \eqref{eq 18b}. The Appendices \ref{s:Sobsp} and \ref{s:Fuc} give the necessary mathematical preparations on the Sobolev spaces and include the main derivations of Fuchsian system used in Sec. \ref{s:case1} and Sec. \ref{s:case3}.


	\section{A non-Fourier-based proof of standard Jeans instability}\label{s:case2}
		In this section, we present a non-Fourier-based proof of standard linearized Jeans instability in the expanding Newtonian universe, based on the technique of \textit{separation of variables}, and in the next section, we develop a method, without calculating the solution, to estimate the behavior of the solution, which can be applied to the nonlinear case. By separating the spatial and time variables, we obtain a Helmholtz equation of the spatial variable and a time-dependent variable coefficient ordinary differential equation (i.e., ODE) with respect to the time. The Helmholtz equation can be solved by the well-known methods and the difficult part is to solve the time-dependent variable coefficient ODE. In the following statement,  in order to simplify the calculations and emphasize the key idea of the proof,
		let us fix  $\gamma= \frac{4}{3}$ in the equation of state \eqref{e:eos}, and  denote $\mathring{\varrho}:=\varrho|_{t=1}$ and $\mathring{\varrho}_\mu:=(\del{\mu}\varrho)|_{t=1}$ ($\mu=0,\cdots,3$).
		
		We assume $\lambda\leq 0$ and $g_\lambda(\boldsymbol{q})$ are eigenvalues and eigenfunctions (which is not equal to $0$ identically) of the Laplace operator $\Delta$ on a manifold\footnote{If we take $M=\Tbb^3$ or various special regions and coordinates, we can calculate $\lambda$ and $g_\lambda(\boldsymbol{q})$; see, for example, \cite{Chavel1984}. } $M$, i.e., $g_\lambda(\boldsymbol{q})$ solves the  Helmholtz equation $\Delta g=\lambda g$ for the eigenvalues $\lambda$ on $M$. In this condition, without loss of generality, we also assume the initial data of density satisfies \footnote{For general data $\mathring{\varrho}=\phi(\boldsymbol{q})$ and $\mathring{\varrho}_0=\psi(\boldsymbol{q})$, we can use the following technique of separation of variables, as usual with the help of the linear superposition principle, to further derive the solution.  }
		\begin{equation}\label{data2}
			\mathring{\varrho}=g_\lambda(\boldsymbol{q})  \AND \mathring{\varrho}_0=\frac{2}{3}g_\lambda(\boldsymbol{q}) 
		\end{equation}
		for a specific eigenvalue
		 $\lambda$ satisfies
		\begin{equation}\label{e:jeanscri0}
			\lambda \in \Bigl[- \frac{25( 6\pi G )^{\frac{1}{3}}}{48\kappa},0\Bigr].
		\end{equation}
		In the next section, we will prove the solution of Eq. \eqref{eq 9} is
		\begin{align}\label{e:rhoexp}
			\varrho(t,\boldsymbol{q})
			=& \Bigl(\frac{1}{2}-\frac{5}{2 \sqrt{36 \lambda\tilde{\kappa}+25}}\Bigr)t^{\frac{2}{3}- \frac{5+\sqrt{25+36\lambda\tilde{\kappa}}}{6} }g( \boldsymbol{q}) +\Bigl(\frac{5}{2 \sqrt{36 \lambda\tilde{\kappa}+25}}+\frac{1}{2}\Bigr)t^{\frac{2}{3}- \frac{5-\sqrt{25+36\lambda\tilde{\kappa}}}{6} }g( \boldsymbol{q}) ,
		\end{align}
	where, we recall, $ \tilde{\kappa}=\frac{4}{3} \kappa(\frac{1}{6\pi G})^{\frac{1}{3}}$.
		Furthermore,  from \eqref{e:rhoexp}, we can see if  $\lambda$ satisfies the \textit{Jeans criterion}:
		\begin{equation}\label{e:jeanscri}
			\lambda \in \Bigl(- \frac{ ( 6\pi G )^{\frac{1}{3}}}{2\kappa},0\Bigr],
		\end{equation}
then
		for every $(t,\boldsymbol{q})\in \Rbb_{>0}\times \mathbb{D}$ where the domain  $\mathbb{D}:=\{\boldsymbol{q}\in \Tbb^3 \;|\; g_\lambda (\boldsymbol{q})  >0\}$, $\varrho$ is growing.
		
 	We point out that \eqref{e:jeanscri} can be viewed as the \textit{Jeans' criterion} in this context. For example, if $M=\Tbb^3$, i.e., for $\boldsymbol{q}\in \Tbb^3$,  $g_\lambda(\boldsymbol{q})=   e^{i\boldsymbol{k}\cdot \boldsymbol{q} }$, then $\lambda=-(k_1^2+k_2^2+k_3^2)$. Then \eqref{e:jeanscri} implies
 	\begin{equation*}
 		k_1^2+k_2^2+k_3^2  < \frac{( 6\pi G )^{\frac{1}{3}}}{2\kappa},
 	\end{equation*}
 	which, in fact, is the \textit{Jeans' criterion}.

	Now let us calculate $\varrho$ is given by  \eqref{e:rhoexp}.
		First we assume the density $\varrho$ has the following separation of variables, i.e., assume there is a function $f(t)$, such that:
		\begin{align}
			\left\{
			\begin{aligned}
				\varrho(t,\boldsymbol{q})=& t^{\frac{2}{3}}f(t) g( \boldsymbol{q}) , \\
				\varrho_{0}(t,\boldsymbol{q})=&f_{0}(t)g( \boldsymbol{q})t^{\frac{2}{3}}+\frac{2}{3}t^{-\frac{1}{3}}f(t)g( \boldsymbol{q}),\\
				\varrho_i(t,\boldsymbol{q})=&t^{\frac{2}{3} }f(t)g_i( \boldsymbol{q}),  \label{eq 12}
			\end{aligned}
			\right.
		\end{align}
		where $g_i:=\del{i}g$.
		Then substituting \eqref{eq 12} into \eqref{eq 9} (noting $\gamma=4/3$), and by using \eqref{data2}, we arrive at a Helmholtz equation $\Delta g=\lambda g$ and a system of $f(t)$,
		\begin{align}\label{e:feq}
			\left\{
			\begin{aligned}
				&\ddot{f}+\frac{8}{3t}\dot{f}-\lambda\tilde{\kappa} \frac{1}{t^2} f=0, \\
				&f|_{t=1}=1 \AND f_0|_{t=1}=0.
			\end{aligned}
			\right.
		\end{align}
	The standard theory on Laplace operator and Helmholtz equations gives the explicit expressions of the eigenvalues $\lambda$ and eigenfunctions $g_\lambda$ in various regions and coordinates, we omit the details, but only focus on the equation of $f(t)$ here.
	
		By letting $f_0:= \del{t}f$, direct calculations imply \eqref{e:feq} can be rewritten as a first-order system,
		\begin{align}
			\left\{
			\begin{aligned}
				&\partial_{t}f_0+\frac{8}{3t}f_{0}-\lambda \tilde{\kappa}  \frac{1}{t^2} f=0,\\
				&\partial_{t}f=f_0,\\
				&f|_{t=1}=1 \AND f_0|_{t=1}=0. \label{eq 16}
			\end{aligned}
			\right.
		\end{align}	
		We introduce a time transform $\tau=\frac{1}{t}$, new  variables $\mathbf{f}_0(\tau):= -t  f_0(t)$ and $\mathbf{f}(\tau):=f(t)$, and denote $V(\tau):=(\mathbf{f}_0(\tau), \mathbf{f}(\tau))^T$.
		Then, we reexpress \eqref{eq 16} into the Fuchsian form,\footnote{This sheds some lights on the methods in Sec. \ref{s:case1}. }
		\begin{equation}
			\partial_{\tau} V=\frac{1}{\tau} \p{
				\frac{5}{3}  & \lambda\tilde{\kappa} \\ 1 & 0} V .  \label{eq 14c}
		\end{equation}
		We denote $\tilde{V}:=(h_1,h_2)$ and let
		\begin{align}\label{e:VV}
			\tilde{V}:=\p{1 & \frac{-5+\sqrt{25+36\lambda\tilde{\kappa}
				}}{6}  \\
				1 & \frac{-5-\sqrt{25+36\lambda\tilde{\kappa}}}{6} }V=\p{\mathbf{f}_0+ \frac{-5+\sqrt{25+36\lambda\tilde{\kappa}
				}}{6}  \mathbf{f} \\ \mathbf{f}_0+ \frac{-5-\sqrt{25+36\lambda\tilde{\kappa}
				}}{6}  \mathbf{f} }.
		\end{align}
		Then interpreting the \eqref{eq 14c} in terms of $\tilde{V}$, we arrive at
		\begin{align}\label{e:eqtilv}
			\del{\tau}\tilde{V}=\frac{1}{\tau} \p{
				\frac{\sqrt{36 \lambda\tilde{\kappa}+25}+5}{6}  & 0 \\
				0 & \frac{5-\sqrt{36 \lambda\tilde{\kappa}+25}}{6}  }\tilde{V}.
		\end{align}
		Noting this system \eqref{e:eqtilv} has decoupled the equations of $h_1$ and $h_2$, we can now solve these ODEs, respectively. Further calculations for equations of $h_1$ and $h_2$, respectively, we obtain
		\begin{equation}\label{e:heq}
			\del{\tau}(\tau^{ \frac{-5-\sqrt{25+36\lambda\tilde{\kappa}}}{6} }h_1)=0 \AND \del{\tau}(\tau^{ \frac{-5+\sqrt{25+36\lambda\tilde{\kappa}}}{6} }h_2)=0.
		\end{equation}
		Solving equations \eqref{e:heq} and replacing, via \eqref{e:VV}, $h_1$ and $h_2$ by $\mathbf{f}_0$ and $\mathbf{f}$, we reach
		\begin{align}
			\mathbf{f}_0+ \frac{-5+\sqrt{25+36\lambda\tilde{\kappa}
			}}{6}  \mathbf{f} =&\frac{-5+\sqrt{25+36\lambda\tilde{\kappa}
			}}{6}t^{- \frac{5+\sqrt{25+36\lambda\tilde{\kappa}}}{6} }, \label{e:ff1} \\
			\mathbf{f}_0+ \frac{-5-\sqrt{25+36\lambda\tilde{\kappa}
			}}{6}  \mathbf{f}=&\frac{-5-\sqrt{25+36\lambda\tilde{\kappa}
			}}{6}t^{- \frac{5-\sqrt{25+36\lambda\tilde{\kappa}\lambda\tilde{\kappa}}}{6} }.  \label{e:ff2}
		\end{align}
		Solving \eqref{e:ff1}--\eqref{e:ff2} and expressing the solutions in terms of $f$ and $f_0$ yield
		\begin{align*}
			f(t)=&\Bigl(\frac{1}{2}-\frac{5}{2 \sqrt{36 \lambda\tilde{\kappa}+25}}\Bigr)t^{- \frac{5+\sqrt{25+36\lambda\tilde{\kappa}}}{6} }+\Bigl(\frac{5}{2 \sqrt{36 \lambda\tilde{\kappa}+25}}+\frac{1}{2}\Bigr)t^{- \frac{5-\sqrt{25+36\lambda\tilde{\kappa}}}{6} },  \\
			f_0(t)=&-\frac{3\lambda\tilde{\kappa}}{\sqrt{25+36\lambda\tilde{\kappa}}} \bigl(t^{- \frac{11+\sqrt{25+36\lambda\tilde{\kappa}}}{6} } -t^{- \frac{11-\sqrt{25+36\lambda\tilde{\kappa}}}{6} } \bigr).
		\end{align*}
		Then with the help of \eqref{eq 12}, we conclude \eqref{e:rhoexp}.


	\section{Fuchsian methods for the linear Jeans instability}\label{s:case1}
	Inspired by the expression of \eqref{eq 14c}, in this section,
	 we develop the Fuchsian formulation associating to the Eq. \eqref{eq 9} and estimate the behavior of the solution to the Eq. \eqref{eq 9} without a detailed solving of this equation. This section presents the main idea of this article and it can be applied to the nonlinear case (see Sec. \ref{s:case3}).  However, before the nonlinear setting, we use the linear case to state the method in order to convey the main idea, but leave the slightly nonlinear Jeans instability to the next section.
	
	Recall that, in the last section, we have obtained the solution to \eqref{eq 9} under the data given by the certain eigenfunctions of the Laplacian [the eigenvalues are restrained by \eqref{e:jeanscri0}]. For a large class of the general initial data, we are able to use the linear superposition principle to construct the solution.  In order to simplify the calculations and highlight the main structures, we, in this section, only focus on $M=\Tbb^3$ and the initial data $\mathring{\varrho}$ perturbing around $g_{0}=\text{Constant}>0$ (i.e., the eigenvalue $\lambda=0$). It is clear $g_0$ satisfies the Jeans criterion \eqref{e:jeanscri}.  In other words, we assume the Fourier series of $\mathring{\varrho}$ is  $\mathring{\varrho}=c_0+\sum_{0\neq \boldsymbol{k}\in \Zbb^n} c_{\boldsymbol{k}}e^{i\boldsymbol{k}\cdot \boldsymbol{q} }$, where $c_{\boldsymbol{k}}$ are the Fourier coefficients, then we require $c_0>0$. We also need $c_{\boldsymbol{k}}$ to be bounded in the certain sense given in the following and $\mathring{\varrho}_\mu$ have similar requirements. It turns out we can obtain the increasing solutions for every $\gamma>1$ (this can also be seen in the previous calculations in Sec. \ref{s:case2}) instead of the only fixed $\gamma=4/3$ if the data is near a positive constant $g_0$. In fact, this initial profile of density with vanishing wave number (i.e., $g_0>0$) leads to the essential accretion.
	
	In order to state the main statement of this case rigorously, we have to introduce a type of function spaces, called \textit{Sobolev spaces}, used widely in the analysis of the nonlinear partial differential equations. We present two definitions in Appendix \ref{s:Sobsp} and one of them is strongly related to the high order derivatives, while the other is from the aspect of Fourier transforms. 
	In order to find a way to prove Jean's instability for the full nonlinear Euler--Poisson system, complex function spaces such as Sobolev spaces are required. 
	It turns out that the conditions and methods given in \textit{Sobolev spaces} and the following \textit{Fuchsian system} can be generated to the slightly nonlinear Jean's instability and are very promising for the full nonlinear Euler--Poisson system.

 Now let us present the main statement of this section:

		Suppose $s\in \Rbb_{\geq 3}$ and $\gamma>1$ are constants and  $\mathring{\varrho}:=\varrho|_{t=1}$ and $\mathring{\varrho}_\mu:=(\del{\mu}\varrho)|_{t=1}$ ($\mu=0,\cdots,3$). Let the initial data of the density satisfy an estimate
		\begin{equation}\label{e:data}
			\Bigl\|\mathring{\varrho}-\frac{\beta}{2}\Bigr\|_{H^s(\Tbb^3)}+\Bigl\|\mathring{\varrho}_0-\frac{\beta}{3}\Bigr\|_{H^s(\Tbb^3)}+\|\mathring{\varrho}_i\|_{H^s(\Tbb^3)}\leq \beta_0 ,
		\end{equation}
		where $0<\beta<+\infty
		$ is any given constant and $\beta_0>0$ is a small enough constant [the detailed range of $\beta_0$ is given by \eqref{e:beta0}]. Then we will prove  the solution of Eq. \eqref{eq 9} satisfies
		\begin{equation*}
		\frac{1}{4} \beta t^{\frac{2}{3}}  \leq \varrho\leq \frac{3}{4} \beta t^{\frac{2}{3}}
		\end{equation*}
		for every $(t,\boldsymbol{q})\in \Rbb_{>0}\times \Tbb^3$.

Now let us prove this statement.
The main idea is  to transform \eqref{eq 9} into the Fuchsian form of \eqref{e:modl} since we have a well-controlled solution of these Fuchsian formulations (see Appendix \ref{s:Fuc}). To achieve this purpose, let us first introduce the following transformation:
		\begin{align}
			\left\{
			\begin{aligned}
				w(t,\boldsymbol{q})&:=\varrho(t,\boldsymbol{q})-\frac{1}{2}\beta t^{\frac{2}{3}},\\
				w_{0}(t,\boldsymbol{q})&:=\partial_{t}w(t,\boldsymbol{q})=\partial_{t}\varrho(t,\boldsymbol{q})-\frac{1}{3}  \beta  t^{-\frac{1}{3}},
				\\
				w_{i}(t,\boldsymbol{q})&:=\partial_{i}w(t,\boldsymbol{q})=\partial_{i}\varrho(t,\boldsymbol{q}) .\label{eq 10}
			\end{aligned}
			\right.
		\end{align}
		Substituting \eqref{eq 10} into \eqref{eq 9}, we reexpress \eqref{eq 9} into a first-order system,
		\begin{align}
			\left\{
			\begin{aligned}
				&\partial_{t}w_0+\frac{4}{3t} w_{0}-\tilde{\kappa}t^{-2\gamma+\frac{2}{3}}\delta^{ij}\partial_{j}w_{i}-\frac{2}{3 t^2} w=0,\\
				&\tilde{\kappa}t^{-2\gamma+\frac{2}{3}}\delta^{ik}\partial_{t}w_{i}-\tilde{\kappa}t^{-2\gamma+\frac{2}{3}}\delta^{ik}\partial_{i}w_{0}=0,\\
				&\partial_{t}w=w_0.\label{eq 11}
			\end{aligned}
			\right.
		\end{align}
		Then rescaling the above set of variables $(w_0,w_i,w)$ and letting $\tau=\frac{1}{t}\in (0,1]$, we denote
		\begin{align}
			\left\{
			\begin{aligned}
				u(\tau,\boldsymbol{q})&:=\frac{\sqrt{6}}{3}t^{-\frac{2}{3}}w(t,\boldsymbol{q})= \frac{\sqrt{6}}{3}t^{-\frac{2}{3}}\varrho(t,\boldsymbol{q})-\frac{\sqrt{6}}{6}  \beta  , \\
				u_{0}(\tau,\boldsymbol{q})&:=t^{\frac{1}{3}}w_0(t,\boldsymbol{q})= t^{\frac{1}{3}}\del{t}\varrho (t,\boldsymbol{q})-\frac{1}{3}  \beta   ,\\
				u_{i}(\tau,\boldsymbol{q})&:=t^{\frac{2}{3}-\gamma}w_i(t,\boldsymbol{q})=t^{\frac{2}{3}-\gamma}\del{i}\varrho(t,\boldsymbol{q})  .\label{eq 12a}
			\end{aligned}
			\right.
		\end{align}
		Replacing the variables \eqref{eq 12a} into \eqref{eq 11}, we arrive at	
		\begin{align}
			\left\{
			\begin{aligned}
				&\partial_{\tau}u_{0}+\tilde{\kappa}\tau^{\gamma-\frac{7}{3}}\delta^{ij}\partial_{j}u_{i}=\frac{1}{\tau}\Bigl(u_0-\frac{\sqrt{6}}{3}u\Bigr),\\
				&\tilde{\kappa}\delta^{ik}\partial_{\tau}u_{i}+\tilde{\kappa}\tau^{\gamma-\frac{7}{3}}\delta^{ik}\partial_{i}u_{0}=\frac{\tilde{\kappa}(\gamma-\frac{2}{3})}{\tau}\delta^{ik}u_i,\\
				&\partial_{\tau}u=\frac{1}{\tau}\Bigl(\frac{2}{3}u-\frac{\sqrt{6}}{3}u_{0}\Bigr).\label{eq 13}
			\end{aligned}
			\right.
		\end{align}
		In terms of the matrix formulations, \eqref{eq 13} becomes
		\begin{equation}
			B^0\partial_{\tau}\mathbf{U}+\tau^{\gamma-\frac{7}{3}}B^i\partial_{i}\mathbf{U}=\frac{1}{\tau}\mathcal{B}\mathbb{P}\mathbf{U},\label{eq 14a}
		\end{equation}
		where $\mathbf{U}:=(u_0,u_j,u)^T$ and $B^0$, $B^i$, $\mathcal{B}$, and $\mathbb{P}$ are constant matrices, i.e. 	
		\begin{gather}
			B^0=\left(                	\begin{matrix}
				1 &  & \\
				& \tilde{\kappa}\delta^{jk} & \\
				&  & 1
			\end{matrix}
			\right),
			\quad
			B^i=\left(                 \begin{matrix}
				0 & \tilde{\kappa}\delta^{ij} & 0\\
				\tilde{\kappa}\delta^{ik} & 0 & 0\\
				0 & 0 & 0
			\end{matrix}
			\right),   \label{e:Bcoef1}
			\\
			\mathcal{B}=\left(
			\begin{matrix}
				\frac{5}{3} &  & \\
				& \tilde{\kappa}(\gamma-\frac{2}{3})\delta^{ik} & \\
				&  & \frac{5}{3}
			\end{matrix}
			\right),
			\quad
			\mathbb{P}=\left(
			\begin{matrix}
				\frac{3}{5} & 0 & -\frac{\sqrt{6}}{5}\\
				0 & \delta^j_i & 0\\
				-\frac{\sqrt{6}}{5} & 0 & \frac{2}{5}
			\end{matrix}
			\right).   \label{e:Bcoef13}
		\end{gather}

		We can, by directly using the model Eq. \eqref{e:modl} and Theorem \ref{t:Fuc} in Appendix \ref{s:Fuc}, Sobolev’s inequality (see Theorem \ref{t:Sobembd} in Appendix \ref{s:Sobsp}), it follows, that the solution $\mathbf{U}=(u_0,u_j,u)^T$ exists on the time interval $\tau\in(0,1]$ and satisfies
		\begin{equation}\label{e:Uest}
			\Vert \mathbf{U} (\tau)\Vert_{L^\infty}\leq  C_s \Vert \mathbf{U} (\tau) \Vert_{H^s}\leq  C_s \Vert \mathring{\mathbf{U}} \Vert_{H^s}
		\end{equation} 	
	where $C_s>0$ is the Sobolev constant from Theorem \ref{t:Sobembd}.
		By \eqref{e:Uest},  the initial data \eqref{e:data} and the transformation \eqref{eq 12a}, we obtain for $\tau\in(0,1]$,
		\begin{equation}\label{e:data2}
			\|\mathbf{U}(\tau)\|_{L^\infty} \leq C_s \|\mathring{\mathbf{U}}\|_{H^s} \leq \frac{\sqrt{6}}{3}\beta_0 C_s \leq \frac{\sqrt{6}}{12}\beta.
		\end{equation}
		It, with the help of \eqref{e:data2} and the transform \eqref{eq 12a},  implies that
		\begin{equation*}
			\frac{1}{4} \beta t^{\frac{2}{3}}  \leq \varrho\leq \frac{3}{4} \beta t^{\frac{2}{3}}.
		\end{equation*}
 This completes the statement and implies $\varrho$ increases in the order $\sim t^{\frac{2}{3}}$.


	\section{Nonlinear analysis of the Jeans instability}\label{s:case3}
	In this section, we consider the Jeans equation with a slightly nonlinear term given by \eqref{eq 18b}, i.e.,
	\begin{equation}
		\ddot{\varrho}+\frac{4}{3t} \dot{\varrho}-\tilde{\kappa}t^{-2\gamma+\frac{2}{3}}\Delta\varrho-\frac{2}{3t^2}\varrho=(\gamma-1)\tilde{\kappa}t^{-2\gamma+\frac{2}{3}}\frac{D^i\varrho D_i\varrho}{1+\varrho}.  \label{eq 18}
	\end{equation}
We point out that the nonlinear term in the right hand of Eq. \eqref{eq 18} is not chosen arbitrarily, we select this nonlinear  term since it appears in the fully nonlinear Jeans equations (derived by the fully nonlinear Euler--Poisson system).

In this section, we are going to prove the similar statement to the above section. That is, suppose $s\in \Rbb_{\geq 3}$ and $\gamma>1$ are constants and  $\mathring{\varrho}:=\varrho|_{t=1}$ and $\mathring{\varrho}_\mu:=(\del{\mu}\varrho)|_{t=1}$ ($\mu=0,\cdots,3$). Let the initial data of the density satisfy an estimate
\begin{equation}\label{e:data1}
	\Bigl\|\mathring{\varrho}-\frac{\beta}{2}\Bigr\|_{H^s(\Tbb^3)}+\Bigl\|\mathring{\varrho}_0-\frac{\beta}{3}\Bigr\|_{H^s(\Tbb^3)}+\|\mathring{\varrho}_i\|_{H^s(\Tbb^3)}\leq \beta_0 ,
\end{equation}
where $0<\beta<+\infty
$ is any given constant and $\beta_0>0$ is a small enough constant [the detailed range of $\beta_0$ is given by \eqref{e:beta0}]. Then we will prove  the solution of Eq. \eqref{eq 18} satisfies
\begin{equation} \label{e:rhoest2}
	\frac{1}{4} \beta t^{\frac{2}{3}}  \leq \varrho\leq \frac{3}{4} \beta t^{\frac{2}{3}}
\end{equation}
for every $(t,\boldsymbol{q})\in \Rbb_{>0}\times \Tbb^3$ as well.

	Similar to Sec. \ref{s:case1}, we substitute \eqref{eq 10} into \eqref{eq 18} and obtain the following first-order system:
	\begin{align}
		\left\{
		\begin{aligned}
			&\partial_{t}w_0+2Hw_{0}-\tilde{\kappa}t^{-2\gamma+\frac{2}{3}}\delta^{ij}\partial_{j}w_{i}-\frac{2}{3 t^2} w=\frac{(\gamma-1)\tilde{\kappa}t^{-2\gamma+\frac{2}{3}}\delta^{ij}w_i w_j}{w+\frac{1}{2} \beta  t^{\frac{2}{3}}+1},\\
			&\tilde{\kappa}t^{-2\gamma+\frac{2}{3}}\delta^{ik}\partial_{t}w_{i}-\tilde{\kappa}t^{-2\gamma+\frac{2}{3}}\delta^{ik}\partial_{i}w_{0}=0,\\
			&\partial_{t}w=w_0.\label{eq 19}
		\end{aligned}
		\right.
	\end{align}
	By similar steps as above, replacing the variables in \eqref{eq 19} by \eqref{eq 12a}, we have the following matrix expression:
	\begin{equation}
		B^0\partial_{\tau}\mathbf{U}+\tau^{\gamma-\frac{7}{3}}B^i\partial_{i}\mathbf{U}=\frac{1}{\tau}\mathcal{B}\mathbb{P}\mathbf{U}+\frac{1}{\tau}H,\label{eq 14b}
	\end{equation}
	where $\mathbf{U}:=(u_0,u_j,u)^T$, $B^0$, $B^i$, $\mathcal{B}$, and $\Pbb$ are constant matrices given by \eqref{e:Bcoef1}--\eqref{e:Bcoef13}, and $H$ is
	\begin{align}\label{e:Bcoef3}
		H=\biggl(-\frac{2\tilde{\kappa}(\gamma-1)\delta^{ij} u_i u_j}{\sqrt{6}u+ \beta  +2\tau^{\frac{2}{3}}}, 0, 0 \biggr)^T  .
	\end{align}  	
This has turned \eqref{eq 18} into the Fuchsian form in \eqref{e:modl}. Similarly using Theorem \ref{t:Fuc}, we conclude the above statement.


	
	\section{Conclusions and discussions}
    We emphasize that we only consider a nonlinear toy model \eqref{eq 18} in Sec. \ref{s:case3} of this article due to the complexity of the fully nonlinear Jeans instability, and, in this toy model, eventually, the nonlinearity does not affect the solution too much, no matter if the initial perturbations are large or small. In fact, this is the exact reason why we study this relatively simple nonlinearity at the first step. In other words, our proof rigorously indicates that the nonlinearity of the nonlinear equation \eqref{eq 18} in  Sec. \ref{s:case3} can not overwhelm the linear effects in the evolutions. Of course, for the fully nonlinear case, we believe the solution will be affected significantly by other nonlinearities of the system and will be studied in the  near future. We emphasize that this paper currently \textit{aims} to rigorously prove this linear dominant phenomenon for the given slightly nonlinear equations. Since we are directly solving and estimating the exact solution  of the nonlinear equation \eqref{eq 18}, \textit{the transition from linear to the nonlinear case in fact has already been included in the detailed proof of the Fuchsian method} (see Appendix \ref{s:Fuc}) of solving this slightly nonlinear equation, as we said before; in fact, this nonlinear equation is dominated by the linear one, and we use a small perturbation of the  increasing mode of the linear Jeans analysis given in  Sec. \ref{s:case2} and  Sec. \ref{s:case1} (by a non-Fourier-based proof) and use the Fuchsian method to conclude an estimate which implies the growing mode of  the solution to this nonlinear equation is exactly the same as the linear one. Roughly speaking, in the proof of the Fuchsian model in Appendix \ref{s:Fuc}, it implies the nonlinearity of \eqref{eq 18} can be fully absorbed by the ``good'' $1/\tilde{\tau}$ singular term (which helps the proof); thus its behavior is fully like the linear one. No matter how large the data is, the rigorous proof yields that the nonlinear solution is indeed dominated by the linear one, but we highlight that this is true only for this special nonlinearity given in \eqref{eq 18} rather than the fully nonlinear Euler--Poisson system. For general nonlinear case, we have to use more delicate method and we are attempting to generalize this Fuchsian method to indicate how the nonlinear behavior dominates.  On the other hand,  our proof is independent of the original classical Jeans' method, thus our initial density perturbations could be large (mathematically allowable) just representing it is already in the procedure of the mass accretions [but note it may not represent the physical reality due to the fact that we are only considering the slightly nonlinear toy model \eqref{eq 18}]. 
    The estimate \eqref{e:rhoest2} bounds the exact solution for all time rather than a small time since the proof of this slightly nonlinear equation \eqref{eq 18} is rigorous without approximations, but it only works for this equation \eqref{eq 18} instead of the fully nonlinear Euler--Poisson equations.
    Therefore, the estimates \eqref{e:rhoest2} of the exact solution still can not conclude the fully nonlinear Jeans instability since we can not include all other nonlinearities of Euler--Poisson equations.

In summary, we have rigorously proved the slightly nonlinear Jeans instability in the expanding Newtonian universe and obtained that the increasing rate of $\varrho$  is not affected by the nonlinearity given by the square of spatial derivatives. However, we conjecture that the increasing rate will be significantly changed by other more difficult nonlinearities and the increasing rate of Jeans instability of the fully nonlinear Euler--Poisson system is very different from $t^{\frac{2}{3}}$ since roughly, speaking, we can see if the density $\varrho$ is increasing, then the coefficient term of the pressure  term $-\frac{c_s^2}{a^2}\Delta\varrho$ in \eqref{eq 7} will change significantly and these effects will feedback to the system leading to faster increasing rates.

Instabilities are very different from the stabilities and thus the stability can be approximated by the solutions of the linearized equations, but the instability of solutions  will deviate from solutions to the linearized equations significantly. Therefore, it is necessary to consider the Jeans instability of the fully nonlinear Euler--Poisson system, and if we can obtain a very fast increasing rate, then we may explain the large inhomogeneities today by small assumptions on the inhomogeneities at the early universe; that is, we may not require strong constraints on the initial spectrum of the perturbations at the early stage of the universe. The proceeding works are on the fully nonlinear Jeans instability both for the Newtonian universe and general relativity.

	\appendix

	\section{Sobolev inequalities and Cauchy-Schwarz inequalities}\label{s:Sobsp}
	In this Appendix, let us very briefly introduce the Sobolev spaces without further mathematical terminologies and details. We refer to, for example,  \cite{Adams2003b,Evans2010} and the references therein for details. We give two equivalent definitions of function spaces $H^{s}(\mathbb{T}^n)$ and explain these two $H^s$ norms which are both used in this article are \textit{equivalent}.  Mathematically, if the norms are equivalent, then it does not matter which norms one uses, they do not affect the analysis since they characterize, in mathematical terminology, the same topology.
	\begin{definition}\label{t:Sobsp}
		When $s$ is a real number, the Sobolev space $H^{s}(\mathbb{T}^n)$ consists of all functions $u$ satisfying
		\begin{equation}\label{e:norm1}
			\Vert u \Vert_{H^s(\mathbb{T}^n)}:=\biggl( \sum_{\xi\in \mathbb{Z}^n} \lvert \hat{u}(\xi) \rvert^{2}(1+\lvert \xi \rvert^{2})^s  \biggr)^{1/2}<\infty,
		\end{equation}
		where $\hat{u}$ denotes the Fourier transform of $u$ and $\Vert u \Vert_{H^s(\Tbb^n)}$ is called the $H^s$-norm of $u$.
	\end{definition}
	Another equivalent definition is given by
	\begin{definition}\label{t:Sobsp2}
		When $s$ is a real number, the Sobolev space $H^{s}(\mathbb{T}^n)$ consists of all functions $u$ satisfying
		\begin{equation}\label{e:norm2}
			\Vert u \Vert_{H^s(\mathbb{T}^n)}:=\biggl( \sum_{\lvert\alpha\rvert\leq s}
			\int_{\mathbb{T}^n}\lvert D^{\alpha} u \rvert^2\ dx\biggr)^{1/2}<\infty,
		\end{equation}
		where $D^{\alpha}u$ is the $\alpha^{\rm{th}}$-weak partial derivative derivative of $u$.
	\end{definition}
	
	$H^s$-norms given by \eqref{e:norm1} and \eqref{e:norm2} are equivalent in the following meaning.
	\begin{definition}\label{t:Sobsp3}
		Suppose $\Vert \cdot \Vert$ and $\Vert \cdot \Vert^{'}$ are two norms on $H^s(\mathbb{T}^n)$. Then we call $\Vert \cdot \Vert$ and $\Vert \cdot \Vert^{'}$ are equivalent if there exist two positive real numbers $C_2\geq C_1>0$ such that
		\begin{equation}\nonumber
			C_1\Vert u \Vert \leq \Vert u \Vert^{'} \leq C_2\Vert u \Vert,\qquad \forall u \in H^s(\mathbb{T}^n).
		\end{equation}
	\end{definition}
	
	\begin{theorem}[Sobolev Embedding Theorem]\label{t:Sobembd}
		If $s>k+n/2$ for some integer $k\geqslant 0$,
		then there is a positive constant $C_s>0$, such that
		\begin{equation}\nonumber
			\sum_{\lvert \alpha \rvert\leqslant k}\Vert \partial^{\alpha}u \Vert_{L^\infty(\Tbb^n)}\leq C_s\Vert u \Vert_{H^s(\Tbb^n)}, \qquad \forall u \in H^s(\Tbb^n),
		\end{equation}
		where $\Vert f\Vert_{L^\infty(\Tbb^n)}:=\sup_{x\in\Tbb^n}|f(x)|$ is the $L^\infty$-norm of $f$ and the constant $C_s$ depends only on $s$, $n$ and $k$.
	\end{theorem}
	
	We present next the following variations of Moser's estimates\footnote{These Moser's estimates can be obtained by theorems in \cite{Taylor2010}, Proposition 3.9, or \cite{Liu2018}, Lemma A.3, with the help of the Sobolev embedding theorem \ref{t:Sobembd}. } which are widely used in the nonlinear PDEs.
	\begin{theorem}[Moser's estimates]\label{t:moser1}
		Let $s\in\Zbb_{\geq 3}$, $F$ be smooth. Then, there are constants $\tilde{C}_m>0$ and $C_m>0$ depending on the maximums of $D^kF$ ($k=1,\cdots s$), such that for $u\in H^s(\Tbb^3)\cap L^\infty(\Tbb^3)$,
		\begin{equation*}
			\|F(u)\|_{H^s} \leq \tilde{C}_m \|u\|_{L^\infty}^{s-1} \|u\|_{H^s} \leq C_m\|u\|_{H^s}^s .
		\end{equation*}
	\end{theorem}
\begin{theorem}[Moser's estimates]\label{t:moser2}
	Suppose $s\in\Zbb_{\geq3}$, $f_i\in H^s(\Tbb^3)$, then there is a constant $C>0$, such that
	\begin{align*}
		\|f_1\cdots f_i\|_{H^s} \leq C\prod^l_{i=1}\|f_i\|_{H^s}.
	\end{align*}
\end{theorem}

The following inequalities are also widely used in the nonlinear PDEs, we list them without proof.
\begin{theorem}[H\"older inequality]\label{t:holder}
	For any functions $f$ and $g$, there is an inequality,
	\begin{equation*}
		|\la f,g\ra| \leq \int_{\Tbb^n}|fg| dx  \leq \Bigl(\int_{\Tbb^n}|f|^2 dx\Bigr)^{\frac{1}{2}}\Bigl(\int_{\Tbb^n}|g|^2 dx\Bigr)^{\frac{1}{2}},
	\end{equation*}
where we denote $\la g_1,g_2\ra:=\int_{\Tbb^3}g_1 g_2 d^3\boldsymbol{q}$.
\end{theorem}
\begin{theorem}[Cauchy–Schwarz inequality]\label{t:cauchy}
	For any vectors $(u_i)$ and $(v_i)$, there is an inequality,
	\begin{equation*}
		  \Bigl(\sum _{i}u_{i}v_{i}\Bigr)^{2}\leq \Bigl(\sum _{i}{u_{i}^{2}}\Bigr)\Bigl(\sum _{i}{v_{i}^{2}}\Bigr).
	\end{equation*}
\end{theorem}


	\section{Fuchsian systems}\label{s:Fuc}
The systems given in \eqref{eq 14a} and \eqref{eq 14b} are in the Fuchsian formulation.
	The initial value problem of quasilinear Fuchsian system has been investigated in a series works by Oliynyk, Beyer, Olvera-Santamar{\'{\i}}a, and the first author of this paper in \cite{Oliynyk2016a,Liu2018,Liu2018b,Beyer2020,LW2021a}.
	In this Appendix, we only focus on the simplest case serving to this article.
	
	Let us, by introducing $\hta=-\tau$, denote \eqref{eq 14a} and \eqref{eq 14b} uniformly by
		\begin{equation}\label{e:modl}
		B^0\partial_{\hta}\mathbf{U}-(-\hta)^{\gamma-\frac{7}{3}}B^i\partial_{i}\mathbf{U}=\frac{1}{\hta}\mathcal{B}\mathbb{P}\mathbf{U}+\frac{\epsilon}{\hta}H, \quad \text{for}\quad \hta\in[-1,0)
	\end{equation}
where $\epsilon=1$ is corresponding to \eqref{eq 14b} and $\epsilon=0$ is corresponding to \eqref{eq 14a}, $\mathbf{U}:=(u_0,u_j,u)^T$, $B^0$, $B^i$, $\mathcal{B}$, and $\Pbb$ are constant matrices given by \eqref{e:Bcoef1}--\eqref{e:Bcoef13}, and $H$ is given by \eqref{e:Bcoef3}.  For this system, we have the following conclusion.

	\begin{theorem}\label{t:Fuc}
		Suppose $s\in\Zbb_{\geq3}$, the initial data  $\mathring{\mathbf{U}}:=\mathbf{U}|_{\hta=-1}\in H^s(\Tbb^3)$, and if there is a small enough [see \eqref{e:beta0}] constant $\beta_0>0$, such that the initial data has an upper bound $\|\mathring{\mathbf{U}}\|_{H^s} \leq \beta_0$,
		then there exists a solution $\mathbf{U}\in\mathnormal C^1([-1,0)\times\mathbb{T}^3)$ to the equation \eqref{e:modl}, such that $\|\mathbf{U}(\hta)\|_{H^s}$ is a  non-increasing function for $\hta\in[-1,0)$ and thus satisfies the estimate
		\begin{equation}\label{e:estFuch}
			\Vert \mathbf{U}(\hta) \Vert_{H^s(\Tbb^3)}
			\leq\Vert \mathring{\mathbf{U}} \Vert_{H^s(\Tbb^3)}
		\end{equation}
		for $\hta\in[-1,0)$.
	\end{theorem}
\begin{remark}
	In fact, one can calculate if taking
	\begin{equation}\label{e:beta0}
		\beta_0\leq \min\biggl\{\frac{ \beta}{8 C_s},\, \frac{1}{2}\Bigl(\frac{\lambda_0}{\epsilon C_m}\Bigr)^{\frac{1}{s+1}}\biggr\} \AND \lambda_0:=\frac{\min\{\frac{5}{3},\gamma-\frac{2}{3}\}}{\max\{1,\frac{1}{\tilde{\kappa}}\}}
	\end{equation}
where $\beta$ is given by \eqref{e:data}, $C_m$ is given by \eqref{e:Hest} from the following proof and $C_s$ is the Sobolev constant given later, then the above theorem holds.
\end{remark}

	\begin{proof}
		The local in time existence and uniqueness  can be obtained, as a special case, by the standard theory of nonlinear hyperbolic equations (for example, see \cite{Taylor2010}, Chap. 16). Now let us obtain the estimates \eqref{e:estFuch} from the Eq. \eqref{e:modl} and derive the solution that exists for $\hta\in[-1,0)$.  Firstly, let us find a time $T\in(-1,0)$ in the local existence interval such that  $\|\mathbf{U}(\hta)\|_{H^s}  \leq 2\beta_0$ for $\hta\in[-1,T)$ (this can be done by the continuity of solution $\mathbf{U}$ and the initial data $\|\mathring{\mathbf{U}}\|_{H^s} \leq \beta_0$).  This leads, by the Sobolev embedding Theorem \ref{t:Sobembd}, to
		\begin{equation}\label{e:usup}
			\|\mathbf{U}(\hta)\|_{L^\infty} \leq C_s\|\mathbf{U}(\hta)\|_{H^s}  \leq 2C_s \beta_0
		\end{equation}
		where $C_s>0$ is the Sobolev constant and it is independent of $\hta$.
		
		Acting on the both sides of \eqref{e:modl} by $D^{\alpha}(B^{0})^{-1}$ yields
		\begin{equation}
			\partial_{\hta} D^{\alpha} \mathbf{U}-(-\hta)^{\gamma-\frac{7}{3}}(B^{0})^{-1}B^{i}\partial_{i} D^{\alpha} \mathbf{U}=\frac{1}{\hta}( B^{0})^{-1}\mathcal{B} D^{\alpha}\mathbb{P} \mathbf{U}+ \frac{\epsilon}{\hta} (B^0)^{-1} D^\alpha H(\hta,\mathbf{U}, \Pbb  \mathbf{U})  . \label{eq A.3}
		\end{equation}
		By denoting $\la g_1,g_2\ra:=\int_{\Tbb^3}g_1 g_2 d^3\boldsymbol{q}$, and using $\la D^\alpha \mathbf{U},\cdot \ra$ to act on \eqref{eq A.3}, with the help of the integration by parts,
		we arrive at\footnote{Note, by \eqref{e:Bcoef1}--\eqref{e:Bcoef13} and direct calculations, $\Pbb^T=\Pbb$, $[(B^0)^{-1}\mathcal{B},\Pbb]=0$ and $\Pbb^2=\Pbb$, then we have $\langle\mathnormal{D^{\alpha}} \mathbf{U},(\mathnormal{B}^{0})^{-1}\mathcal{B}\mathnormal{D^{\alpha}}\mathbb{P}\mathbf{U} \rangle =\langle\mathnormal{D^{\alpha}}\mathbb{P}\mathbf{U},(\mathnormal{B}^{0})^{-1}\mathcal{B}\mathnormal{D^{\alpha}}\mathbb{P}\mathbf{U} \rangle $. }
		\begin{align}			   &\frac{1}{2}\partial_{\hta}\langle\mathnormal{D^{\alpha}}\mathbf{U},\mathnormal{D^{\alpha}}\mathbf{U}\rangle -(-\hta)^{\gamma-\frac{7}{3}}\langle\mathnormal{D^{\alpha}}\mathbf{U},(\mathnormal{B}^{0})^{-1}\mathnormal{B}^{i}\partial_{i}\mathnormal{D^{\alpha}}\mathbf{U}\rangle\notag  \\
		&\hspace{1cm}= \frac{1}{\hta}\langle\mathnormal{D^{\alpha}}\mathbb{P}\mathbf{U},(\mathnormal{B}^{0})^{-1}\mathcal{B}\mathnormal{D^{\alpha}}\mathbb{P}\mathbf{U} \rangle  + \frac{\epsilon}{\hta} \la  D^\alpha \mathbf{U}, (B^0)^{-1} D^\alpha H(\hta,\mathbf{U},\Pbb \mathbf{U}) \ra .\label{e:engeq}
		\end{align}
Note
\begin{equation*}
	2\langle\mathnormal{D^{\alpha}}\mathbf{U},(\mathnormal{B}^{0})^{-1}\mathnormal{B}^{i}\partial_{i}\mathnormal{D^{\alpha}}\mathbf{U}\rangle=\int_{\Tbb^3}\partial_{i}\bigl((\mathnormal{D^{\alpha}}\mathbf{U})^T (\mathnormal{B}^{0})^{-1}\mathnormal{B}^{i}\mathnormal{D^{\alpha}}\mathbf{U}\bigr) d^3\boldsymbol{q}=0,
\end{equation*}
due to the divergence theorem and the integration by parts.
We denote, by recalling Definition \ref{t:Sobsp2}, the energy norm
		\begin{equation}\nonumber
			\Vert \mathbf{U} \Vert_{H^s}^{2}=\sum_{\lvert\alpha\rvert\leq  s}\langle\mathnormal{D^{\alpha}}\mathbf{U},\mathnormal{D^{\alpha}}\mathbf{U}\rangle.
		\end{equation}
Let us, by Theorems \ref{t:moser1} and \ref{t:moser2}, and noting $\tau\in(-T,1]$ (i.e., $\hta\in[-1,T)$), firstly estimate $H$ [we recall $H$ is given by \eqref{e:Bcoef3}],
\begin{align}\label{e:Hest}
	\|H\|_{H^s} =&\biggl\|-\frac{2\tilde{\kappa}(\gamma-1)\delta^{ij} u_i u_j}{\sqrt{6}u+\beta+2\tau^{\frac{2}{3}}} \biggr\|_{H^s} \notag  \\
	\leq& C\biggl\|-\frac{2\tilde{\kappa}(\gamma-1) }{\sqrt{6}u+ \beta  +2\tau^{\frac{2}{3}}} \biggr\|_{H^s}\|\Pbb \mathbf{U}\|^2_{H^s}\notag  \\
	\leq& C_m \|u\|^s_{H^s}\|\Pbb \mathbf{U}\|^2_{H^s}\leq   C_m (2\beta_0)^s \|\Pbb \mathbf{U}\|^2_{H^s}
\end{align}	
where the constant $C_m$ can be chosen (large enough) to be a time-independent constant for the finite $\tau \in(-T,1]$.
Note that in the above calculations, we have used, by \eqref{e:usup} and \eqref{e:beta0},  $\sqrt{6}u+\beta\geq -\sqrt{6}\|u\|_{L^\infty}+\beta>-2\sqrt{6}C_s\beta_0+\beta\geq \bigl(1-\frac{\sqrt{6}}{4}\bigr)\beta>0$.
	
Then let us, by using \eqref{e:Hest}, H\"older and Cauchy–Schwarz inequalities (see Theorem \ref{t:holder} and \ref{t:cauchy}), estimate the last term of \eqref{e:engeq} for $\hta\in[-1,T)$,
		\begin{align}\label{e:Hest1}
			&\Bigl|\sum_{|\alpha|\leq s} \la  D^\alpha \mathbf{U}, (B^0)^{-1} D^\alpha H(\hta,\mathbf{U},\Pbb \mathbf{U})	 \ra \Bigr|  \notag  \\
		    \overset{\text{H\"older inequality}}{\leq } & \max\Bigl\{1,\frac{1}{\tilde{\kappa}}\Bigr\} \sum_{|\alpha|\leq s}  \biggl(\int_{\Tbb^3}|D^\alpha\mathbf{U}|^2d^3\boldsymbol{q}\biggr)^{\frac{1}{2}} \biggl(\int_{\Tbb^3}|D^\alpha H|^2d^3\boldsymbol{q} \biggr)^{\frac{1}{2}}\notag  \\
			\overset{\text{Cauchy inequality}}{\leq } &  \max\Bigl\{1,\frac{1}{\tilde{\kappa}}\Bigr\} \biggl(\sum_{|\alpha|\leq s} \int_{\Tbb^3}|D^\alpha\mathbf{U}|^2d^3\boldsymbol{q}\biggr)^{\frac{1}{2}}  \biggl(\sum_{|\alpha|\leq s} \int_{\Tbb^3}|D^\alpha H|^2d^3\boldsymbol{q}\biggr)^{\frac{1}{2}} \notag  \\
			=&  \max\Bigl\{1,\frac{1}{\tilde{\kappa}} \Bigr\} \| \mathbf{U}\|_{H^s}  \|  H\|_{H^s} \notag  \\
			\leq   &  C_m (2\beta_0)^s  \max\Bigl\{1,\frac{1}{\tilde{\kappa}}\Bigr\}\|\mathbf{U}\|_{H^s} \|\Pbb \mathbf{U}\|_{H^s}^2 \notag  \\
			\leq   &  C_m (2\beta_0)^{s+1} \max\Bigl\{1,\frac{1}{\tilde{\kappa}}\Bigr\}  \|\Pbb \mathbf{U}\|_{H^s}^2 .
		\end{align}
In addition, we note
\begin{equation}\label{e:singtm}
	\la D^\alpha \Pbb\mathbf{U},(B^0)^{-1}\mathcal{B} D^\alpha\Pbb\mathbf{U}\ra\geq \min\Bigl\{\frac{5}{3}, \Bigl(\gamma-\frac{2}{3}\Bigr)\Bigr\}\|\Pbb\mathbf{U}\|^2_{H^s} .
\end{equation}

		Next,  by using $2\sum_{\lvert\alpha\rvert\leq  s} $ acting on \eqref{e:engeq} and taking $\beta_0$ small enough [at least, such that \eqref{e:beta0} holds], with the help of \eqref{e:Hest1} and \eqref{e:singtm}, we reach the estimate (note $\hta<0$ in the following calculations)
		\begin{equation}\nonumber
			\partial_{\hta}\Vert \mathbf{U}(\hta) \Vert_{H^s}^{2}\leq \frac{2}{\hta}\max\Bigl\{1,\frac{1}{\tilde{\kappa}}\Bigr\}\Bigl(\lambda_0-\epsilon C_m (2\beta_0)^{s+1} \Bigr)\Vert  \mathbb{P}\mathbf{U} \Vert_{H^s}^{2}\leq  0.
		\end{equation}
	Therefore, the function $\|\mathbf{U}(\hta)\|_{H^s}$ is a nonincreasing function with respect to time $\hta\in[-1,T)$.
From this, we find that
	for any $\hta\in[-1,T)$,
	\begin{align}\label{e:Uest1}
		\Vert \mathbf{U}(\hta) \Vert_{H^s} \leq 	\Vert \mathring{\mathbf{U}}\Vert_{H^s}  \leq  \beta_0.
	\end{align}

By continuous extensions, i.e., letting $\mathbf{U}(T):=\lim_{\hta\rightarrow T}\mathbf{U}(\hta)$, we extend the solution $\mathbf{U}(\hta)$ for  $\hta\in[-1,T)$ to $\hta=T$ satisfying $\|\mathbf{U}(T)\|_{H^s}\leq \beta_0$. Then using $\mathbf{U}(T)$ as the initial data, we are able to apply above derivations to the time interval $[T,-1+2(T+1))$, and conclude $\mathbf{U}$ exists for $\bigl[-1,-1+2(T+1)\bigr)\times \Tbb^3$ and $\|\mathbf{U}(\hta)\|_{H^s}$ is a nonincreasing function in this time interval. Furthermore, repeatedly using the above method and calculations, it yields that the solution can be extended to the time interval $\hta\in\bigl[-1,-1+\ell(T+1)\bigr)$ (for $\ell\in\Zbb_{>0}$) and the estimate \eqref{e:Uest1} holds for $\hta\in\bigl[-1,-1+\ell(T+1)\bigr)$.
Then there must be a finite step $\ell=\ell_0$ such that $\ell_0$ is the smallest number satisfying $-1+\ell_0(T+1)\geq 0$. This means\footnote{That is, at the last $\ell_0$-th step, we proceed above derivations in $\hta\in \bigl[-1+(\ell_0-1)(T+1),0\bigr)$. } we can enlarge the time interval of the existence of the solution to $\hta\in[-1,0)$ and $\|\mathbf{U}(\hta)\|_{H^s}$  is a  nonincreasing function for $\hta\in[-1,0)$ (which implies \eqref{e:estFuch} as well). Then we complete the proof.
	\end{proof}

	\section*{Acknowledgement}
	This work is partially supported by the Fundamental Research Funds for the Central Universities, HUST: $5003011036$, $5003011047$.

	\bibliographystyle{unsrt}
	\bibliography{Reference_Chao}

\end{document}